# NEW EXOTIC MINIMAL SETS FROM PSEUDO-SUSPENSIONS OF CANTOR SYSTEMS

J. P. BOROŃSKI, ALEX CLARK, AND P. OPROCHA

ABSTRACT. We develop a technique, pseudo-suspension, that applies to invariant sets of homeomorphisms of a class of annulus homeomorphisms we describe, Handel-Anosov-Katok (HAK) homeomorphisms, that generalize the homeomorphism first described by Handel. Given a HAK homeomorphism and a homeomorphism of the Cantor set, the pseudo-suspension yields a homeomorphism of a new space that admits a homeomorphism that combines features of both of the original homeomorphisms. This allows us to answer a well known open question by providing examples of hereditarily indecomposable continua that admit homeomorphisms of intermediate complexity. Additionally, we show that such examples occur as minimal sets of volume preserving smooth diffeomorphisms of 4-dimensional manifolds. We also use our techniques to exhibit the first examples of minimal, uniformly rigid and weakly mixing homeomorphisms in dimension 1, and these can also be realized as invariant sets of smooth diffeomorphisms of a 4-manifold. Until now the only known examples of spaces that admit minimal, uniformly rigid and weakly mixing homeomorphisms were modifications of those given by Glasner and Maon in dimension at least 2.

## 1. INTRODUCTION

This paper is motivated by recent advances in surface dynamics, where peculiar invariant sets and attractors play a prominent role. Many of these efforts focus on understanding dynamical systems on cofrontiers, or other invariant continua, and the study of the interplay between the complicated topology of such sets and the complexity of the associated dynamics. There is a very long, and well documented history of this line of study, starting with the work of Poincaré. Subsequent milestones include the discovery of strange attractors with rotational chaos by Birkhoff [12], the study of fixed points in acyclic invariant continua by Cartwright and Littlewood [32], Smale's and Williams' work on expanding hyperbolic attractors and the construction of DA-attractors [74, 76, 78] and Barge and Diamond's work on the topology and dynamics of inverse limits of graphs [5]. The research on this topic has intensified in recent years in rotation theory of surfaces, which has significantly advanced our understanding of the phenomenon frequently referred to as strange or exotic attractors and minimal sets [56]. Moreover, the work of T. Jäger and his collaborators suggests that the dynamics and topology of peculiar continua might be intertwined with the unsolved case of Franks-Misiurewicz Conjecture; see e.g. [8, 50, 51]. In [23] and [24] a systematic study has been carried out of attractors obtained from embeddings of inverse limits. The main method of those papers extends upon the Barge-Martin-Brown dynamical embedding theorem [16, 6]. We propose a new method of constructing 1-dimensional attractors by suspending cofrontier dynamics onto Cantor set minimal systems. As a model for an exotic cofrontier we choose R.H. Bing's pseudo-circle, as topologically it is the most pathological object in the class, being hereditarily indecomposable. Recall that a metric space with at least two points is a *continuum* if it is compact and connected, and a *cofrontier* $X$ is a planar continuum whose complement is the union of two connected open sets for each of which $X$ is the boundary. A continuum $X$ is *indecomposable* if it is not the union of two of its proper subcontinua, and when additionally every subcontinuum of a continuum $X$ is indecomposable, then $X$ is *hereditarily indecomposable*.

However, as we make clear, the method applies in a much more general setting. As one of the applications of our method, we investigate the dynamics of hereditarily indecomposable continua that appear with intermediate complexity. We say that a map $T\colon X \to X$ of a metric space $X$ is of

---

*intermediate complexity* when its topological entropy satisfies $0 < h_{top}(T) < +\infty$. There are many examples of hereditarily indecomposable continua occurring as invariant sets of smooth dynamical systems with zero entropy. Handel [44] and later Herman [47] using the so–called Anosov-Katok fast approximation method provided examples of smooth maps of two dimensional manifolds with Bing's pseudo-circle [11] as a minimal set. In these examples the maps can be chosen to be $C^\infty$– smooth, volume preserving diffeomorphisms. Handel's example has the additional property that a slight perturbation yields the existence of a $C^\infty$–smooth diffeomorphism of the annulus for which the pseudo-circle is an attracting minimal set [44] with the same homeomorphism when restricted to the invariant pseudo-circle. Kennedy and Yorke [52, 53, 54] described $C^\infty$–smooth dynamical systems in dimensions greater than 2 with uncountably many minimal pseudo-circles and satisfying the stability condition that any small $C^1$ perturbation also has this same property. More recently hereditarily indecomposable continua made their appearences in complex dynamics. In [29] Chéritat extended Herman's construction to show that the pseudo-circle occurs as the boundary of a Siegel disk for a holomorphic map. Subsequently, in the context of Eremenko's conjecture, Rempe constructed examples of entire functions whose Julia sets have components that are pseudo-arcs [68], thereby strengthening the results from [73].

Many of these examples indicate that hereditarily indecomposable continua occur quite naturally with zero entropy within smooth systems. This leads to the following initial well known question that one can find posed in [62, 19].

**QUESTION 1.1.** *Is there a hereditarily indecomposable continuum $X$ that admits a homeomorphism of intermediate complexity?*

A variant of that question was raised by Barge in 1989, who asked whether the pseudo-arc (the arc-like hereditarily indecomposable continuum) admits homeomorphisms with arbitrary value of entropy (see [61]). These questions have important implications for the possible ways that hereditarily indecomposable continua occur as invariant sets of smooth systems. Ito [49] has shown that any $C^1$ diffeomorphism on a compact finite-dimensional Riemannian manifold must have finite topological entropy. Although the questions have remained unsolved until now, there were some partial results indicating that the answer to both may be in the negative. This was proved first for shifts of inverse limits on the arc by Mouron [63], and then on general topological graphs $G$ by Boroński and Oprocha [19]. Namely, if the hereditarily indecomposable continuum $X$ is presented as an inverse limit of a single map $f\colon G \to G$, then the corresponding shift map has either 0 or infinite topological entropy. Furthermore, if the graph $G$ is a circle, then the only admissible value of entropy is infinity. Since attractors can often be presented as such shifts (see [4, 6, 23, 24, 76]), these results indicate that there might be an obstruction to the existence of hereditarily indecomposable continua as invariant sets of intermediate complexity. And by the result of Ito, this might lead one to conjecture that any hereditarily indecomposable continuum occurring as an invariant set of a smooth map would be of 0 entropy as in the previously known examples. This would also have some consequences for potential approaches to the remaining unsolved case of Franks-Misiurewicz Conjecture (see [8]).

However, here we establish the following result (its proof is presented in Section 8).

**THEOREM 1.2.** *There exist an hereditarily indecomposable continuum $\Psi_C$ satisfying:*

(1) *$\Psi_C$ occurs as an invariant minimal set with intermediate complexity for a volume preserving diffeomorphism $F$ of a 4 dimensional manifold and*
(2) *The restriction $F|\Psi_C$ is weakly mixing.*

We construct these examples using the technique of pseudo-suspension that we introduce in Section 3. This technique can be considered as an adaptation of the classic suspension (or mapping cylinder) construction as it applies to homeomorphisms of the annulus. In particular, for each transitive homeomorphism $h\colon C \to C$ of the Cantor set $C$, we construct a one–dimensional continuum that supports self homeomorphisms that can be considered as the lifts of homeomorphisms of the pseudo-circle or more general cofrontiers. When the homeomorphism $h\colon C \to C$ is minimal and the cofrontier is the pseudo-circle, the constructed continuum is hereditarily indecomposable.

There is another important observation related to our construction. It would seem that since Handel's homeomorphism is obtained as a limit of rigid rotations that its dynamics would also be similar to a rotation. In particular, it is natural to expect that it would be uniformly rigid and distal. While the first claim is true, the second one is false. Following observation of Thurston, presented by Handel in [44] we easily obtain that Handel's example is topologically weakly mixing (see Corollary 5.4). This way we obtain an example of minimal, weakly mixing and uniformly rigid homeomorphism of the pseudo-circle, and by our method we are also able to extend it to other one–dimensional hereditarily indecomposable continua. It is worth mentioning that the only known construction of weakly mixing, minimal and uniformly rigid homeomorphism originated from [41] and is a consequence of Baire category techniques of construction of minimal skew-products developed by Glasner and Weiss in [42]. The technique from [41] immediately provides examples on tori (of any dimension $n \geq 2$) and was further generalized to obtain examples on other surfaces (e.g. see [80]). However, previously there was no known example in dimension one.

Relating to the above, one of the starting points to this research was a question raised by Joe Auslander after authors' talk on the results in [18], at the AIMS conference in Madrid in 2014, on the existence of proximal pairs in Handel homeomorphism. Initially it seemed almost evident that no such pairs may exist and all pairs are distal, but it seemed of interest to determine whether an adaptation of Anosov-Katok fast approximation method, employed by Handel, could yield new examples with such pairs or Li-Yorke chaos. It was during such considerations that the authors were led to realize that Handel's homeomorphism is weakly mixing and that this is intertwined with the lack of semi-conjugacy to a circle rotation. We also obtain results that establish a link between the nature of the rotation set of a self homeomorphism of the pseudo-circle and its topological entropy. In particular, we have the following (see Section 10 for a proof):

**THEOREM 1.3.** *Let $\Psi \subset \mathbb{A}$ be an essential pseudo-circle attracting all the points from $\mathrm{Int}\,\mathbb{A}$ and assume that $H \colon \mathbb{A} \to \mathbb{A}$ is a homemomorphism with a nondegenerate rotation set. Then $h_{top}(H|_\Psi) = h_{top}(H) = +\infty$.*

This extends a recent result of Passeggi, Potrie and Sambarino [67], who showed that a nondegenerate rotation interval of an annulus homeomorphism $H$ on an invariant attracting cofrontier $K$ implies positive entropy of $H|_K$.

2. PRELIMINARIES

We assume that the reader is familiar with basic notions of topological dynamics, such as minimality, topological mixing and topological entropy (for more details we refer the reader to monographs [40, 60, 77]).

2.1. **Handel-Anosov-Katok homeomorphisms.** In [44] Handel employed the so-called Anosov–Katok fast approximation method to obtain a cofrontier in the annulus $\mathbb{A}$ with special properties. In this section we will present minimal ingredients that are essential to obtain examples sharing key properties with those obtained by Handel.

2.1.1. *Notation for the Annulus.* The universal covering of $\mathbb{A} = [0,1] \times \mathbb{R}/\mathbb{Z}$ is given by $q \colon [0,1] \times \mathbb{R} \to \mathbb{A}$, where $q = \mathrm{id}_{[0,1]} \times p$ and $p$ is the quotient map $\mathbb{R} \to \mathbb{R}/\mathbb{Z}$. We endow $[0,1] \times \mathbb{R}$ with the metric given by adding the Euclidean metrics in the factors, and $\mathbb{A}$ is given the rotation invariant metric $d$ for which $q$ is a local isometry. Here and in what follows $\pi_i$ denotes the projection of $\mathbb{A}$ onto its $i$-th factor ($i = 1, 2$). We use $\rho$ to denote the standard complete $C^0$ metric in the space of all self homeomorphisms of $\mathbb{A}$.

We say that a homeomorphism $H \colon \mathbb{A} \to \mathbb{A}$ is a *Handel-Anosov-Katok homeomorphism*, or a *HAK homeomorphism* for short, if $H$ is a uniform limit (relative to $\rho$) of a sequence of annulus homeomorphisms $(H_n)_{n \in \mathbb{N}}$ satisfying the following approximation scheme.

We assume first that we are given a decreasing sequence $A_n$ ($n \in \mathbb{N}$) of closed annuli satisfying for each $n \in \mathbb{N}$, $A_{n+1} \subset A_n$ and $A_1 \subset \operatorname{Int} \mathbb{A}$. Moreover, we presume given a convergent series $\sum_{i=1}^{\infty} \varepsilon_i < +\infty$, $(\varepsilon_i > 0)$ and denote $\gamma_n = \sum_{i=n}^{\infty} \varepsilon_i$.

(1) For each $n \in \mathbb{N}$ there exist a homeomorphism $f_n \colon A_n \to \mathbb{A}$ with small vertical fibers; that is, the fibers of the projection $\pi_1 \circ f_n$ have diameters less than $\varepsilon_n/2$.
(2) For each $n \in \mathbb{N}$ we have a rational rotation $R_n \colon \mathbb{A} \to \mathbb{A}$ such that for $\delta_n$ corresponding to $\varepsilon_n/4$ in the uniform continuity of $f_n^{-1}$ there exists a rectangle $D_n = [0,1] \times [0, \alpha_n] \subset \mathbb{A}$ such that the sets $R_n^i(D_n), R_n^j(D_n)$ either have disjoint interiors or are equal, $\bigcup_{i \in \mathbb{Z}} R_n^i(D_n) = \mathbb{A}$ and $\alpha_n < \delta_n$. We denote the period of $R_n$ by $p_n$; that is, $R_n^{p_n} = \operatorname{id}$ and $p_n$ is the smallest positive integer satisfying this property. We also require that $p_{n+1} \geq p_n$.
(3) For each $n \in \mathbb{N}$ there is a homeomorphism $g_n \colon \mathbb{A} \to \mathbb{A}$ such that $\rho(g_n, \operatorname{id}) < \varepsilon_n$ and $g_n|_{A_n} = f_n^{-1} \circ R_n \circ f_n$ and $g_n|_{\mathbb{A} \setminus A_{n-1}} = \operatorname{id}$ where we put $A_0 = \mathbb{A}$,
(4) $H_n = g_n \circ g_{n-1} \circ \ldots \circ g_1$,
(5) $H_n(A_{n+1}) = A_{n+1}$ and
(6) $\rho(H_n^i, H_{n+1}^i) < \varepsilon_n$ for $i = 1, \ldots, q_n$ and $q_n = m p_n$ for some $m \geq 1$.

If we want our examples to be smooth (say $C^{\infty}$), it suffices to replace $\rho$ by the analogous $C^n$-metric $\rho_n$ in (4) and (6) for each $n$. It follows from (1)–(6) that $\rho(H_n, H_{n+j}) \leq \gamma_n$ and so $(H_n)_{n \in \mathbb{N}}$ is a Cauchy sequence, implying $H_n$ converges to a $C^k$ diffeomorphism $H$ in the respective space.

First observe that $H_n(A_n) = A_n$ because clearly $g_1(A_1) = A_1$ and by (5) and induction we have
$$H_n(\mathbb{A} \setminus A_n) = g_n(H_{n-1}(\mathbb{A} \setminus A_n)) = g_n(\mathbb{A} \setminus A_n) = \mathbb{A} \setminus A_n$$
which then gives
$$H_{n+m}(\mathbb{A} \setminus A_n) = g_{n+m} \circ \ldots g_{n+1} \circ H_n(\mathbb{A} \setminus A_n) = g_{n+m} \circ \ldots g_{n+1}(\mathbb{A} \setminus A_n) = \operatorname{id}(\mathbb{A} \setminus A_n) = \mathbb{A} \setminus A_n.$$
This immediately gives $H(A_n) = A_n$ for each $n$.

Put $D_n^0 = f_n^{-1}(D_n)$ and for each $i = 1, \ldots, p_n - 1$ denote $D_n^i = g_n^i(D_n^0) = f_n^{-1}(R_n^i(D_n))$. Now for any $z \in D_n^i$ denote by $I_z$ the interval of the form $[0,1] \times \{a\}$ such that $f_n(z) \in I_z \subset R_n^i(D_n)$. Since $R_n$ is an isometry, for any $x, y \in D_n^i$ there are $p \in I_x$ and $q \in I_y$ such that $d(p,q) \leq \operatorname{diam} \pi_1(D_n) < \delta_n$. But then by (1) we have:
$$\begin{aligned} d(x,y) &\leq \operatorname{diam}(f_n^{-1}(I_x)) + \operatorname{diam}(f_n^{-1}(I_y)) + d(f_n^{-1}(p), f_n^{-1}(q)) \\ &< \varepsilon_n/4 + \varepsilon_n/4 + \varepsilon_n/2 = \varepsilon_n. \end{aligned}$$
This implies that

(7) for each $n$ and each $0 \leq j < p_n$ we have $\operatorname{diam} D_n^j < \varepsilon_n$.

Denote
$$\Psi = \bigcap_{n=1}^{\infty} A_n.$$
Clearly $\Psi$ is closed, connected and $H(\Psi) = \Psi$. It is easy to see that $\Psi$ is a cofrontier, and it can always be made an attracting set for $H$ using the same technique as Handel [44] by replacing the condition in (4) by $H_n(A_{n+1}) \subset \operatorname{Int} A_{n+1}$. For this reason we shall call $\Psi$ a *HAK attractor* of the HAK homeomorphism $H$.

**DEFINITION 2.1.** *A homeomorphism $H \colon X \to X$ of the metric space $X$ is* uniformly rigid *if there is a subsequence of iterates $H^{n_i}$ ($n_i \nearrow \infty$) such that the sequence $H^{n_i}$ converges uniformly to $\operatorname{id}_X$ as $i \to \infty$.*

If we fix any $x \in \Psi$, then for each $n$ there is $i$ such that $x \in D_n^i$. For simplicity of notation assume that $x \in D_n^0$. For any $y \in \Psi$ there is $0 \leq j < p_n$ such that $y \in D_n^j$ and then $H_n^j(x) \in D_n^j$ which gives by (6) and (7) that
$$d(H^j(x), y) \leq d(H^j(x), H_n^j(x)) + d(H_n^j(x), y) \leq \gamma_n + \varepsilon_n.$$

This shows that orbit of every $x \in \Psi$ is dense in $\Psi$. Furthermore, since $H_n^{p_n} = \text{id}$ we see that $d(H^{p_n}(x), x) < \gamma_n$. This shows that $H$ is uniformly rigid. Additionally observe that by (6) we have

(8) if $x \in D_n^j$ then $H^i(x) \in B_{\gamma_n}(D_n^{(i+j) \bmod p_n})$ for $j = 0, \ldots, q_n$, where $B_r(S)$ denotes the set of all points of $\mathbb{A}$ within $r$ of $S$ relative to $d$.

By the above discussion, we have the following.

**PROPOSITION 2.2.** *If $H$ is a HAK homeomorphism satisfying* (1)–(6), *then $H$ is minimal and uniformly rigid.*

2.2. **Uniform rigidity and rotation.** Let $\widetilde{\pi}_1$ denote projection onto the first coordinate of $[0,1] \times \mathbb{R}$. For any subset $S \subset \mathbb{A}$, we let $\widetilde{S}$ denote $q^{-1}(S)$. For simplicity, for $(t, r) \in [0,1] \times \mathbb{R}$ and $j \in \mathbb{Z}$ we will simply write $(t, r) + j$ to denote $(t, r + j)$.

Here we see the consequences that a homeomorphism $H$ of $\Psi \subset \mathbb{A}$ is uniformly rigid, regardless of the origin of $H$.

**LEMMA 2.3.** *If $(\Psi, H)$ is uniformly rigid and $\widetilde{\Psi}$ is connected, then for every $\varepsilon > 0$ and every $N$ there are $k > 0$ and $j \in \mathbb{Z}$ such that for each $y \in \widetilde{\Psi}$ we have*
$$|\widetilde{\pi}_1(\tilde{H}^k(y)) - \widetilde{\pi}_1(y) - j| < \varepsilon.$$

*Proof.* Without loss of generality we may assume that $\varepsilon < 1/4$. By uniform rigidity, there is $k > 0$ such that $|H^k(x) - x| < \varepsilon$ for each $x \in \Psi$. Since the universal cover $q \colon \widetilde{\Psi} \to \Psi$ is a local isometry, for each $y \in \widetilde{\Psi}$ there is $j_y \in \mathbb{Z}$ such that
$$\left|\widetilde{\pi}_1(\tilde{H}^k(y)) - \widetilde{\pi}_1(y) - j_y\right| < \varepsilon.$$

Since $\varepsilon < 1/4$, the vector $j_y$ is uniquely determined and the map $g \colon y \mapsto j_y$ is continuous. But as $\widetilde{\Psi}$ is connected, $g$ must be constant. □

**THEOREM 2.4.** *If $(\Psi, H)$ is uniformly rigid and $\widetilde{\Psi}$ is connected, then for every $\varepsilon > 0$ there are $N > 0$ and $\alpha_\varepsilon \in \mathbb{R}$ such that for every $n \geq N$ and every $y \in \widetilde{\Psi}$ we have*
$$\left|\frac{\widetilde{\pi}_1(\tilde{H}^n(y)) - \widetilde{\pi}_1(y)}{n} - \alpha_\varepsilon\right| < \varepsilon.$$

*Proof.* Let $j$ and $k$ correspond to $\varepsilon/4$ as provided by Lemma 2.3. For every $y$ in $\widetilde{\Psi}$ and every positive integer $i$ we have
$$\begin{aligned}|\widetilde{\pi}_1(\tilde{H}^{ik}(y)) - \widetilde{\pi}_1(y) - ij| &\leq |\widetilde{\pi}_1(\tilde{H}^{ik}(y)) - \widetilde{\pi}_1(\tilde{H}^{(i-1)k}(y)) - j| + \ldots + |\widetilde{\pi}_1(\tilde{H}^k(y)) - \widetilde{\pi}_1(y) - j| \\ &< i\varepsilon/4.\end{aligned}$$

Fix any $N > 0$ such that for every $n > N$ and every $i$ satisfying $ik \leq n < (i+1)k$ we have:

- for each $s = 0, 1, \ldots, k$ we have $|\widetilde{\pi}_1(\tilde{H}^s(x)) - \widetilde{\pi}_1(x)| < n\varepsilon/4$ for every $x \in \widetilde{\Psi}$
- $j/N < \varepsilon/4$

Put $\alpha_\varepsilon = \frac{j}{k}$. If we fix any $n > N$ and pick $i$ such that $ik \leq n < (i+1)k$ then
$$\begin{aligned}\left|\frac{\widetilde{\pi}_1(\tilde{H}^n(y)) - \widetilde{\pi}_1(y)}{n} - \alpha_\varepsilon\right| &\leq \left|\frac{\widetilde{\pi}_1(\tilde{H}^n(y)) - \widetilde{\pi}_1(\tilde{H}^{ik}(y))}{n}\right| + \frac{1}{n}\left|\widetilde{\pi}_1(\tilde{H}^{ik}(y)) - \widetilde{\pi}_1(y) - ij\right| + \left|\frac{ij}{n} - \alpha_\varepsilon\right| \\ &< \frac{\varepsilon}{4} + \frac{i\varepsilon}{4n} + \left|\frac{ikj}{nk} - \frac{jn}{nk}\right| \leq \frac{\varepsilon}{2} + \frac{j}{n} \\ &< \varepsilon.\end{aligned}$$

The proof is completed. □

As an immediate consequence of Theorem 2.4 we obtain the following.

**COROLLARY 2.5.** *If $(\Psi, H)$ is uniformly rigid and $\tilde{\Psi}$ is connected then there is well defined rotation number $\alpha$ for all $x \in \Psi$.*

## 3. Pseudo–suspensions

3.1. **Basic construction.** A standard, yet extremely useful technique to construct examples in dynamics is given by suspensions of Cantor dynamical systems. From another point of view, we can view this technique as modification of rigid rotation of the circle to some richer dynamical system that incorporates the dynamics of the Cantor dynamical system. Unfortunately, there are two disadvantages of this approach. First of all, the space supporting the new homeomorphism is relatively simple (a kind of solenoid) and furthermore there are strong restrictions on the dynamics of the resulting systems.it has rotation as a factor.

In this section we will develop a technique which allows us to replace the circle in the suspension construction by a continuum supporting a HAK homeomorphism. That way we will be able to construct spaces with much richer structure. Additionally, the rotation factor is replaced by the HAK homeomorphism as a factor, which leaves much more freedom in finding dynamical properties.

Formally, assume that we are given a cofrontier $\Psi$ as an essentially embedded subspace of the annulus $\mathbb{A}$ and that we are given a minimal homeomorphism $h \colon C \to C$ of the Cantor set $C$. Recall that a continuum $\Psi \subset \mathrm{Int}\,\mathbb{A}$ is *essential* if $\mathbb{A} \setminus \Psi$ has two connected components $\mathcal{U}_+, \mathcal{U}_-$, each of which contains a different component of the boundary of $\mathbb{A}$.

In general we have the principal $\mathbb{Z}$–bundle $\xi = ([0,1] \times \mathbb{R}, q, \mathbb{A})$, where $q : [0,1] \times \mathbb{R} \to \mathbb{A}$ is the universal covering map as before. We then give $C$ the $\mathbb{Z}$–space structure induced by $h$, $n.c \mapsto h^{-n}(c)$, giving us the fiber bundle $\xi[C]$, a different bundle for each choice of $h$. The total space $([0,1] \times \mathbb{R})_C$ of $\xi[C]$ is the quotient space $([0,1] \times \mathbb{R}) \times C/\approx$, where $((s,r),c) \approx ((s',r'),c')$ if and only if $s = s'$ and there is an $n \in \mathbb{Z}$ satisfying $r' = r + n$ and $c' = h^{-n}(c)$, and we denote the quotient map as $Q : ([0,1] \times \mathbb{R}) \times C \to ([0,1] \times \mathbb{R})_C$.

Now we have the following,

$$[0,1] \times \mathbb{R} \times C \xrightarrow{Q} ([0,1] \times \mathbb{R})_C \xrightarrow{q_C} \mathbb{A}$$

where the map $q_C$ is the bundle projection given by $q_C[((s,r),c)] = q((s,r))$, which is well–defined since any representative of the same $\approx$ class differs from $((s,r),c)$ in the first pair of coordinates by a deck transformation of the covering $q$.

**DEFINITION 3.1.** *The subspace $\Psi_C := (q_C)^{-1}(\Psi)$ associated to the homeomorphism $h \colon C \to C$ is the* pseudo–suspension *of $h$.*

The pseudo–suspension of $h$ is the total space of the restricted bundle $\eta := (\Psi_C, q_C, \Psi)$ of the fiber bundle $\xi[C]$. Later in Subsection 3.3 we shall see how $\eta$ can equivalently be regarded as an induced bundle.

3.2. **Lifting homeomorphisms.** Here we describe how to lift a HAK homeomorphism of $\Psi$ to $\Psi_C$ (e.g. Handel's homeomorphism from [44] on the pseudo-circle or the example of Jäger-Koropecki [51]). With the appropriate choice of $h$, the lifted homeomorphism can be chosen to have positive but finite entropy. Handel's homeomorphism has the important additional property that it extends to a homeomorphism $H : \mathbb{A} \to \mathbb{A}$ that is smooth and area preserving.

As the bundle $\xi[C]$ is locally trivial, the notion of bundle here coincides with that in Spanier [75], and so the bundle projection $q_C : ([0,1] \times \mathbb{R})_C \to \mathbb{A}$ is a fibration, meaning it has the homotopy lifting property with respect to all spaces [75, Cor. 14; p. 96]. The maps $q_C$ and $H \circ q_C$ are homotopic maps $([0,1] \times \mathbb{R})_C \to \mathbb{A}$ since $H$ is homotopic to the identity. But since the identity map to which $H$ is homotopic can be lifted, there is a map $H_C$ making the following diagram commute

(3.1)
$$\begin{CD} ([0,1] \times \mathbb{R})_C @>H_C>> ([0,1] \times \mathbb{R})_C \\ @Vq_CVV @VVq_CV \\ \mathbb{A} @>H>> \mathbb{A}. \end{CD}$$

Now $H^{-1}$ is also homotopic to the identity and so can be lifted to a map $H_C^{-1}$ in the same way, and due to the functoriality of the lifting we see that $H_C$ is a homeomorphism. Also, since $H(\Psi) = \Psi$, we obtain that the restriction of $H_C$ is a homeomorphism on $\Psi_C$.

The quotient map $Q$ can be regarded as the map onto the orbit space of the $\mathbb{Z}$ action $\alpha$ on $[0,1]\times\mathbb{R}\times C$ given by
$$n.((t,r),c) = \big((t, n+r), h^{-n}(c)\big).$$

As $\alpha$ acts properly discontinuously, $Q$ is a covering projection [75, Thm. 7, p. 87] in the sense that each point of $([0,1] \times \mathbb{R})_C$ has an open neighborhood $U$ which pulls back by $Q$ to a disjoint union of open subsets, each of which $Q$ maps homeomorphically onto $U$ (that is, $U$ is *evenly covered*). So, again $Q$ has the homotopy lifting property. Thus we obtain again a lifting making the following diagram commute:

(3.2)
$$\begin{CD} [0,1] \times \mathbb{R} \times C @>\widetilde{H_C}>> [0,1] \times \mathbb{R} \times C \\ @VQVV @VVQV \\ ([0,1] \times \mathbb{R})_C @>H_C>> ([0,1] \times \mathbb{R})_C. \end{CD}$$

Notice that the whole homotopy lifts (although the homotopy cannot leave $\Psi_C$ or its lift fixed). Thus, just as the identity map, $\widetilde{H_C}$ does not permute the path components of $[0,1] \times \mathbb{R} \times C$, which are given by $[0,1] \times \mathbb{R} \times \{c\}$.

**DEFINITION 3.2.** *A subset of the pseudo–suspension $\Psi_C$ of the form $Q\left([0,1] \times \mathbb{R} \times \{c\}\right) \cap \Psi_C$ for some $c \in C$ is a* pseudo–component.

Thus we see that the pseudo–components are not permuted by $H_C$.

3.3. **The topology of the pseudo–suspension.** We start by exploring the relation between the pseudo–suspension of a homeomorphism $h : C \to C$ and the suspension of $h$. The universal covering of $\mathbb{R}/\mathbb{Z}$ yields a principal $\mathbb{Z}$–bundle $\zeta = (\mathbb{R}, p, \mathbb{R}/\mathbb{Z})$ where $p$ is the quotient map as before. We then give $C$ the $\mathbb{Z}$–space structure induced by $h$ as before, yielding the fiber bundle $\zeta[C]$. The total space $\mathbb{R}_C$ of $\zeta[C]$ is the suspension of $h$ and is given by $\mathbb{R} \times C/\approx$, where $(r,c) \approx (r',c')$ if and only if there is an $n \in \mathbb{Z}$ satisfying $r' = r + n$ and $c' = h^{-n}(c)$. If we denote the quotient map $Q_h : \mathbb{R} \times C \to \mathbb{R}_C$, then as before we have

$$\mathbb{R} \times C \xrightarrow{Q_h} \mathbb{R}_C \xrightarrow{q_h} \mathbb{R}/\mathbb{Z}$$

where $q_h$ is a bundle projection of $\zeta[C]$. For a fixed $h$, observe then that we have the bundle morphism from $\eta$ to the bundle $\zeta[C]$ determined by the commutative diagram

(3.3)
$$\begin{CD} \Psi_C @>F>> \mathbb{R}_C \\ @Vq_CVV @VVq_hV \\ \Psi @>f>> \mathbb{R}/\mathbb{Z}, \end{CD}$$

where $f: \Psi \to \mathbb{R}/\mathbb{Z}$ is given by projection onto the circle factor, $(s, [r]) \mapsto [r]$ and similarly $F$ maps the class of $((s,r), c)$ to the class of $(r, c)$. In fact, we can alternately regard $(\Psi_C, q_C, \Psi)$ to be the bundle induced from $\zeta[C]$ by the map $f$, [48, Prop. 2.5.5].

The space $\mathbb{R}_C$ is locally homeomorphic to $(0, 1) \times C$ and so is a *matchbox manifold.*

In what follows, by a *proper* subcontinuum of $X$ we mean a subcontinuum $K \subsetneq X$ that contains more than one point. For completeness, we provide a proof of the following well known result.

**PROPOSITION 3.3.** *The only proper subcontinua of the suspension $\mathbb{R}_C$ are arcs.*

*Proof.* Suppose that $K$ is a proper subcontinuum of $\mathbb{R}_C$. Then the complement $\mathbb{R}_C \setminus K$ contains a basic open set $U \approx G \times (-1, 1)$, where $G$ is a clopen subset of $C$. As $G$ is totally disconnected, the first return map $r$ to the section $G \times \{0\}$ of the minimal suspension flow $\phi$ on $\mathbb{R}_C$ is continuous. As $h$ is minimal, $r$ is a homeomorphism. Thus $\mathbb{R}_C$ is homeomorphic to the suspension of $r$ via a homeomorphism $T: \mathbb{R}_C \to \mathrm{Susp}(r)$ that sends $K$ to a homeomorphic copy contained in the neighbourhood $\mathrm{Susp}(r) \setminus h(U)$ that is homeomorphic to $G \times [0, 1]$. But clearly the only continua with more than a point in $G \times [0, 1]$ are arcs. □

We recall the following definition for clarity. We refer the unfamiliar reader to [64, XI].

**DEFINITION 3.4.** *The* composant *of the point $x$ of the space $X$ is the union of all proper subcontinua of $X$ containing $x$.*

In an indecomposable continuum the composants partition the space, each composant is dense and there are uncountably many composants. Proposition 3.3 then has the following known corollary.

**COROLLARY 3.5.** *The path components of $\mathbb{R}_C$ coincide with the composants of $\mathbb{R}_C$.*

**THEOREM 3.6.** *The pseudo–suspension $\Psi_C$ of a minimal homeomorphism $h: C \to C$ is a continuum. If additionally $\Psi$ is the pseudo-circle then the only proper subcontinua of $\Psi_C$ are pseudo-arcs and consequently, $\Psi_C$ is hereditarily indecomposable.*

*Proof.* To show that $\Psi_C$ is compact, we regard it as the bundle induced by the map $f$ in Equation 3.3. As such, $\Psi_C$ is a closed subset of the compact space $\Psi \times \mathbb{R}_C$ ([48, Ch2.5]) and hence compact. To see that $\Psi_C$ is connected, we begin by relating the universal cover of $\Psi$ to our previous constructions. The principal bundle $\xi = ([0, 1] \times \mathbb{R}, q, \mathbb{A})$ corresponding to the universal cover of $\mathbb{A}$ has the restricted bundle $\psi = (\widetilde{\Psi}, q, \Psi)$ where $\widetilde{\Psi} := q^{-1}(\Psi)$. The total space $\widetilde{\Psi}$ is a connected space that admits a two point compactification that is homeomorphic to the pseudo-arc $P$ (e.g. see [58], [9]). Now we can supplement the commutative Diagram 3.3 as follows:

$$
(3.4) \qquad \begin{array}{ccc} \widetilde{\Psi} \times C & \xrightarrow{\widetilde{F}} & \mathbb{R} \times C \\ {\scriptstyle Q} \downarrow & & \downarrow {\scriptstyle Q_h} \\ \Psi_C & \xrightarrow{F} & \mathbb{R}_C. \end{array}
$$

Thus, each pseudo–component of $\Psi_C$ is the continuous bijective image of $\widetilde{\Psi}$ and each pseudo–component is mapped bijectively by $F$ onto a path component of $\mathbb{R}_C$. By the density of the orbit of each $c \in C$ under the action of $h$, we see that each psueudo–component is dense in $\Psi_C$, which is therefore connected.

Now assume that $\Psi$ is the pseudo-circle and let $K$ be a proper subcontinuum of $\Psi_C$. Then $F(K)$ will be a subcontinuum of $\mathbb{R}_C$, which by Proposition 3.3 is an arc $A$ or a point. We shall assume without loss of generality that it is an arc as the other case follows by the same argument. This arc $A$ lifts via $Q_h$ to an arc $\widetilde{A} \subset \mathbb{R} \times \{c\}$ for some $c \in C$. (There are many such lifts, we choose just one.) By construction, $\widetilde{K} := \widetilde{F}^{-1}(\widetilde{A}) \subset \widetilde{\Psi} \times \{c\}$. Moreover, $Q$ maps $\widetilde{K}$ onto $K$ in $\Psi_C$ by the commutativity of the diagram. Now the restriction of $Q$ to $\widetilde{\Psi} \times \{c\}$ is a continuous bijection that is not closed, which is thus a homeomorphisms when restricted to a compact subset. Thus, we finally obtain that $K$ is homeomorphic to a subcontinuum of the pseudo-arc $P$ we obtain by forming the two point compactification of $\widetilde{\Psi} \times \{c\}$. Thus, as any proper subcontinuum of $P$ is homeomorphic to $P$, we obtain that $K$ is homeomorphic to $P$. As $P$ is indecomposable, we obtain that $\mathbb{R}_C$ is hereditarily indecomposable. □

**REMARK 3.7.** *One should bear in mind that when $\Psi$ is the pseudo-circle then each pseudo-component of $\Psi_C$ contains uncountably many distinct composants.*

## 4. Minimality

One important property of HAK homeomorphism $H$ is that it can be uniformly approximated as closely as desired by a periodic homeomorphism for all iterations up to and including the period of the periodic homeomorphism (that is, the identity map). As we have seen in previous section, it then follows directly that $H$ is uniformly rigid.

Thus, for every $\varepsilon > 0$ we can find positive integers $s$ and $k$ such that $\widetilde{H}^s(t,r) \in B((t,r)+k,\varepsilon)$. We shall require a stronger version of this that is satisfied if we increase the convergence rate of the Cauchy sequence $H_n \to H$ of the periodic homeomorphisms $H_n$ obtained by HAK conditions (1)–(6). For example, if we put $q_n = ap_n$ for some integer $a$ in (6) then for any given $\varepsilon > \varepsilon_n$ $p$, we may find positive integers $s$ and $k$ such that $\widetilde{H}^{sj}(t,r) \in B((t,r)+kj,\varepsilon)$ for $j = 1,\ldots,p$. Recall that a homeomorphism $h$ is *totally transitive (minimal)* if for every $k \in \mathbb{Z}$, $h^k$ is transitive (minimal). By standard Baire arguments, if $h \colon C \to C$ is totally transitive, then there exists a point $c \in C$ such that $c$ has dense orbit under $h^k$ for every $k \in \mathbb{Z}$. By the above comments, the assumptions of the following theorem are satisfied if $h$ is totally transitive and values of sequence $q_n$ in HAK conditions are increasing sufficiently fast to ensure appropriate rates of convergence of the approximating homeomorphisms $f_n$.

**THEOREM 4.1.** *Suppose that there are $c \in C$ and $(t,r) \in \hat{\Psi}$ such that for every $\varepsilon > 0$ there are $k, s$ and $p$ such that:*

  *(1) for every $x \in C$ there is $0 \leq i \leq p$ such that $h^{-ik}(c) \in B(x,\varepsilon)$,*
  *(2) for every $0 \leq j \leq p$ we have $\widetilde{H}^{sj}(t,r) \in B((t,r)+kj,\varepsilon)$.*

*Then $H_C$ is transitive and $z = Q((t,r),c) \in \Psi_C$ has a dense forward orbit under the action of $H_C$.*

*Proof.* Fix any $\varepsilon > 0$, any $x \in C$ and let $k, p, s$ be provided by assumptions and let $i$ be provided by (1). Then
$$\widetilde{H}_C^{sj}((t,r),c) \in B((t,r)+kj,\varepsilon) \times \{c\} \approx B((t,r),\varepsilon) \times \{h^{-kj}(c)\} \subset B((t,r),\varepsilon) \times B(x,\varepsilon).$$
This shows that the $\omega$-limit set of $z = Q((t,r),c)$ contains
$$\omega(z, H_C) \supset Q((t,r) \times C) = q_c^{-1}(t,r).$$
But the orbit of $(t,r)$ under $H$ is dense in $\Psi$, which shows that $Q(\{\widetilde{H}^i(t,r) : i = 0,1,\ldots\} \times C)$ is dense in $\Psi_C$ and hence $\Psi_C = \omega(z, H_C)$ since any $\omega$-limit set is closed. This proves that $z$ has dense forward orbit under $H_C$ and hence $H_C$ is transitive. □

It should be noted that in the above, for any fixed estimate we make we are only considering $H$ to be determined up to some finite, uniformly approximating stage $f_n$. This allows us to make meet the requirements for $p$ that might increase without control as $\varepsilon \to 0$. However, this could restrict the rotation numbers allowed for the $H$ we use to those well approximated by the rational numbers given by the rotation numbers of $f_n$.

**COROLLARY 4.2.** *If $h \colon C \to C$ is minimal and weakly mixing and $H_n \to H$ sufficiently fast, then $H_C$ is minimal.*

*Proof.* It is not hard to check that when $h$ is weakly mixing then $(C, h^k)$ is minimal for every $k = 1, 2, \ldots$ (e.g. see [34, 8, p.129]). Then fixing $\varepsilon > 0$ and $k$ we can take the same integers $s, p$ for every point $c \in C$ in (1) in Theorem 4.1. In other words every $z = Q((t, r), c) \in \Psi_C$ has dense forward orbit under the action of $H_C$ which ends the proof. $\square$

## 5. Weak mixing

Here we examine the mixing properties of the homeomorphisms that we construct. Recall that a homeomorphism $h\colon X \to X$ is (topologically) *weak mixing* if for any non–empty open sets $U_i, V_i$ ($i = 1, 2$) there is a $k \geq 1$ such that $h^k(U_i) \cap V_i \neq \emptyset$ ($i = 1, 2$); that is, the system $h \times h : X \times X \to X \times X$ is topologically transitive, in analogy with the notion of weak mixing from ergodic theory. It is known that for a *minimal* homeomorphism of a compact Hausdorff space $h\colon X \to X$ that weak mixing is equivalent to having no non–trivial equicontinuous factors, (see e.g, [34, V, Thm. 1.19, p. 408]). This has a particularly digestible interpretation in terms of eigenfunctions. Recall that a *(continuous) eigenfunction* of $h\colon X \to X$ is a continuous map $\chi : X \to \mathbb{R}/\mathbb{Z}$ satisfying the property that for some $\alpha \in \mathbb{R}/\mathbb{Z}$ with associated rotation $R_\alpha$, $\chi \circ f = R_\alpha \circ \chi$, and in this setting $\alpha$ is the *eigenvalue associated to* $\chi$. Thus, provided that $\chi$ is onto, $\chi$ provides a semi–conjugacy with $R_\alpha$. We then have the following well known theorem, see, e.g., [34, Thm. II 4.19, p. 82].

**THEOREM 5.1.** *Let $h\colon X \to X$ be a minimal homeomorphism of a compact Hausdorff space $X$. Then every eigenfunction of $h$ is constant iff $h$ is weakly mixing.*

It often happens that homeomorphisms of disconnected spaces admit eigenfunctions that are not constant and yet not surjective. Although this is likely well known, for completeness we show that this is not possible in our context.

**LEMMA 5.2.** *If $\chi$ is an eigenfunction of a minimal homeomorphism $h\colon X \to X$ of a compact and connected Hausdorff space $X$ whose image contains more than one point, then $\chi$ is surjective.*

*Proof.* We first observe that by the minimality of $h$ that if $\chi(X)$ contains at least two points, the associated eigenvalue $\alpha$ cannot be the identity element. Suppose then that under our hypotheses that $\chi(X)$ contains two distinct points, that the associated eigenvalue $\alpha$ is not the identity element, but that $\chi$ is not surjective. Then $\chi(X) = A$ for some arc $A$. But since $\chi \circ h = R_\alpha \circ \chi$, $\chi(X)$ must contain $R_\alpha(A)$, which is a distinct arc from $A$ but of the same length, a contradiction. Thus $\chi$ must be surjective. $\square$

**COROLLARY 5.3.** *Let $h\colon X \to X$ be a minimal homeomorphism of a compact and connected Hausdorff space $X$. Then $h$ is not semi–conjugate to a rotation of $\mathbb{R}/\mathbb{Z}$ iff $h$ is weakly mixing.*

Following a suggestion of Thurston, in [44] Handel shows that a homeomorphism $H : \Psi \to \Psi$ of the type he constructs is not semi–conjugate to a rotation of $\mathbb{R}/\mathbb{Z}$ (or any homeomorphism of $\mathbb{R}/\mathbb{Z}$).

**COROLLARY 5.4.** *A Handel homeomorphism $H : \Psi \to \Psi$ is weakly mixing.* $\square$

**REMARK 5.5.** *By Proposition 6.7 in [41], a minimal uniformly rigid homeomorphism of a compact zero dimensional space is always equicontinuous. So any example of a minimal, uniformly rigid and weakly mixing homeomorphism must be at least one–dimensional. In that sense Handel's example is optimal.*

**THEOREM 5.6.** *Let $H_C : \Psi_C \to \Psi_C$ be the lift of a weakly mixing HAK homeomorphism of the cofrontier $\Psi$ to the pseudo–suspension of the weakly mixing homeomorphism $h\colon C \to C$. Then $H_C$ is weakly mixing and in particular $H_C$ does not admit a semi–conjugacy to a rotation of $\mathbb{R}/\mathbb{Z}$.*

*Proof.* By the results of [3], to establish weak mixing, it is sufficient to show that for any nonempty open sets $U, V \subset \Psi_C$ there is an $l > 0$ such that $H_C^l(U) \cap U \neq \emptyset$ and $H_C^l(U) \cap V \neq \emptyset$. We consider first the lift $\tilde{H}$ of $H$ to the universal cover $[0, 1] \times \mathbb{R}$ of $\mathbb{A}$, and take $K$ such that if $\tilde{H}(t, p) = (t', p')$ then $|t' - t| \leq K$, which we know to exist since $\tilde{H}$ is $\mathbb{Z}$ periodic in the second factor. For any nonempty open set $U$ of $[0,1] \times \mathbb{R}$, $U + \mathbb{Z}$ intersects $(0,1)^2$ in a non–empty open set. Fix then any nonempty open sets $U_1, U_2 \subset (0,1)^2$. By the weak mixing of $H$, there are $s > 0$, a nonempty open $U_3 \subset U_1$ and $(t,r) \in (0,1)^2$ such that $(t,r) \in U_1$, $\tilde{H}^s(t,r) \in U_1 + \mathbb{Z}$ and $\tilde{H}^s(U_3) \subset U_2 + \mathbb{Z}$. There exists $a \in \mathbb{Z}$ such that $\tilde{H}^s(t,r) - a \in [0,1]$. Note that $|a| \leq Ks$. Since $H$ is minimal, there exists a syndetic set $A \subset \mathbb{N}$ such that $\tilde{H}^n(t,r) \in U_3 + \mathbb{Z}$ for every $n \in A$. Let $N$ be such that for every $n \in A$ there is $0 < i < N$ such that $n + i \in A$. We may assume that $N > s$. Fix any nonempty open sets $V_1, V_2 \subset C$. By the weak mixing of $h$ there exist points $c_{-3KN}, \ldots, c_{3KN}, c'_{-3KN}, \ldots, c'_{3KN} \in V_1$ and $m > 0$ such that $h^{-m-i}(c_i) \in V_1$ and $h^{-m-i}(c'_i) \in V_2$ for every $|i| \leq 3KN$, see, e.g., [34, II,4.12, p.79]. There exists $l > 0$ such that if we denote $\tilde{H}^l(t,r) = (u,v)$ then $(u,v) \in U_3 + \mathbb{Z}$ and there is $0 \leq i < KN$ such that $m + i \leq v < m + i + 1$. Furthermore, if we put $\tilde{H}^s(u,v) = (u',v') \in U_2 + \mathbb{Z}$ then $|u - u'| \leq sK < KN$. Let $j$ be such that $m + i + j \leq v' < m + i + j + 1$. Taking all the above constants into account we obtain that

$$(\tilde{H}^l(t,r), c_{m+i}) \approx ((u-m-i,v), h^{-m-i}(c_{m+i})) \in U_3 \times V_1 \subset U_1 \times V_1$$
$$(\tilde{H}^l(\tilde{H}^s(t,r)-a), c'_{m+i+j-a}) \approx ((u'-m-i-j, v'), h^{-m-i-j+a}(c'_{m+i+j-a})) \in U_2 \times V_2.$$

Since $Q_C$ is an open map from the product space, this suffices to above mentioned condition to establish weak mixing [3]. □

A natural question to ask at this point is whether the pseudo-circle is the only HAK attractor that admits weak mixing?

**REMARK 5.7.** *There exists a HAK attractor topologically distinct from the pseudo-circle that admits weakly mixing homeomorphisms.*

A detailed description will be given in a separate publication.

## 6. Computation of entropy

We present estimates of topological entropy of homeomorphisms induced on the space $([0,1] \times \mathbb{R})_C$. While in our constructions $H$ is a very special minimal map on the pseudo-circle, we in fact only use that $h_{\text{top}}(H) = 0$ and rotation number is well defined and positive for every point of $\Psi$.

We denote by $r_n(\varepsilon, T)$ the smallest cardinality of $(n, \varepsilon)$-spanning set for $T$ and by $s_n(\varepsilon, T)$ the maximal cardinality of an $(\varepsilon, n)$-separated set for $T$.

**THEOREM 6.1.** *Assume that $H$ is a HAK homeomorphism with a nonzero rotation number $\alpha$ on $\Psi$. Then $h_{top}(H_C) = |\alpha| h_{top}(h)$.*

*Proof.* Replacing $H$ with $H^{-1}$ if necessary, without loss of generality we may assume that $\alpha > 0$. We assume that the distance between points in $[0,1] \times \mathbb{R} \times C$ is the maximum of the distances calculated for the coordinates in $[0,1] \times \mathbb{R}$ and $C$. For simplicity, the metric is always denoted by $d$.

Since $\alpha > 0$ for every sufficiently large positive integer $k$ there is and integer $p$ such that $\frac{p-1}{k} < \alpha < \frac{p+1}{k}$. By Theorem 2.4 we see that there exists $N$ such that for every $n > N$ we have

$$\frac{p-1}{k} < \frac{\pi_1(H^n(x)) - \pi_1(x)}{n} \leq \frac{p+1}{k}$$

for every $x$, and consequently for every $x$ we have

(6.1) $$n(p-1) < \pi_1(H^{kn}(x)) - \pi_1(x) \leq n(p+1).$$

Denote $D = [0,1] \times [-1,2] \times C$ and that (e.g. by (6.1)) we may assume that $k$ is sufficiently large to ensure that

(6.2) $$\tilde{H}^k_C([0,1] \times [-1,\infty) \times C) \subset [0,1] \times [2,\infty) \times C.$$

Fix any $\varepsilon > 0$. Since each sufficiently small set in $([0,1] \times \mathbb{R})_C$ is evenly covered, there exists an open cover $\{U_i\}_{i=1}^s$ of $D$ such that $Q|_{U_i}$ is a homeomorphism onto its image and $\mathrm{diam}(Q(U_i)) < \varepsilon$. Let $\delta$ be a Lebesgue number of the cover $\{U_i\}_{i=1}^s$ (for $D$).

Fix any positive integer $n > N$ and let:

- $A_n$ be the lift to $[0,1]^2$ of an $(n, \delta)$-spanning set of $H^k$ with cardinality $r_n(\delta, H)$,
- $B_n$ be an $(n(p+2), \delta)$-spanning set for $h^{-1}$ with cardinality $r_{n(p+2)}(\delta, h^{-1})$.

Without loss of generality we may assume that the quotient map $q \colon [0,1]^2 \to [0,1] \times (\mathbb{R}/\mathbb{Z}) = \mathbb{A}$ is a local isometry and so for every $y \in [0,1]^2 \cap \tilde{\Psi}$ there is $x \in A_n$ such that $d(\tilde{H}^{ik}(x), \tilde{H}^{ik}(y)) < \delta$ for $i = 0, \ldots, n$. Take any $y \in \Psi_C$ and take its lift
$$((t,u), c) \in ([0,1)^2 \times C) \cap Q^{-1}(y) \subset \tilde{\Psi} \times C$$
Let $(t', u') \in A_n$ and $c' \in B_n$ be such that
$$d(\tilde{H}^{ik}(t,u), \tilde{H}^{ik}(t', u')) < \delta$$
for $i = 0, \ldots, n$ and
$$d(h^{-j}(c), h^{-j}(c')) < \delta$$
for $j = 0, \ldots, r(p+2) - 1$. Note that by (6.1) and (6.2) we have $\tilde{H}^{ki}(t,u) \in [0,1] \times [0, n(p+1) + 2)$ for each $i = 0, \ldots, n$ and therefore there is $0 \leq j < n(p+1) + 2$ such that such that $\tilde{H}^{ik}_C((t,u), c) \subset [0,1] \times [j, j+1) \times C$. Note that $n > N$ is large, hence we may assume that

(6.3) $$j + 2 \leq n(p+1) + 4 < n(p+2).$$

But then $\tilde{H}^{ki}_C((t', u'), c') \subset [0,1] \times [j-1, j+2) \times C$ and so we have
$$\tilde{H}^{ik}_C((t,u), c) \approx (\tilde{H}^{ik}(t,u) - j, h^{-j}(c)) \in D,$$
$$\tilde{H}^{ik}_C((t', u'), c') \approx (\tilde{H}^{ik}(t', u') - j, h^{-j}(c')) \in D.$$

Note that by (6.3) we obtain
$$d(\tilde{H}^{ik}(t,u) - j, h^{-j}(c), \tilde{H}^{ik}(t', u') - j, h^{-j}(c')) \leq \max\{d(\tilde{H}^{ik}(t,u), \tilde{H}^{ik}(t',u')), d(h^{-j}(c), h^{-j}(c'))\}$$
$$\leq \delta$$

which yields $d(H^{ik}_C(y), H^{ik}_C(z)) < \varepsilon$. We obtain that for some sufficiently small $\gamma < \delta$
$$\limsup_{n\to\infty} \frac{1}{n} \log r_n(H^k_C, \varepsilon) \leq \limsup_{n\to\infty} \frac{1}{n} \log(r_n(\delta, H^k) + r_{n(p+2)}(\delta, h^{-1}))$$
$$\leq \limsup_{n\to\infty} \frac{1}{n} \log(r_n(\delta, H^k) + \limsup_{n\to\infty} \frac{1}{n} \log(r_n(\gamma, h^{-p-2})$$

When we pass with $\varepsilon \to 0$ then also $\delta \to 0$ and $\gamma \to 0$, which in the limit gives
$$h_{\mathrm{top}}(H_C) = \frac{1}{k} h_{\mathrm{top}}(H^k_C) \leq h_{\mathrm{top}}(H) + \frac{1}{k} h_{\mathrm{top}}((h^{-1})^{p+2}) = \frac{p+2}{k} h_{\mathrm{top}}(h^{-1}).$$

Increasing $k$ we see that $p/k \to \alpha$ which finally gives
$$h_{\mathrm{top}}(H_C) \leq \alpha h_{\mathrm{top}}(h^{-1}).$$

For the proof of converse implication fix any $\varepsilon > 0$. Let $k$ and $p$ be as before, in particular (6.1) is satisfied. Let $\gamma < \varepsilon$ be such that if $r, t \in C$ and $d(r,t) < \gamma$ then $d(h^i(r), h^i(t)) < \varepsilon$ for $i = -2s, \ldots, 2s$ where $s$ is such that $\tilde{H}^k(D) \subset [0,1] \times [0, s]$. Since each sufficiently small set in $([0,1] \times \mathbb{R})_C$ is evenly covered, there exists $\delta > 0$ such that if $(x, r), (x, t) \in [0,1]^2 \times C$ and $d(Q(x, r), Q(x, t)) < \delta$ then $d(r, t) < \gamma$.

Fix any $n > N$ and let $C_n$ be an $(n(p+1), \varepsilon)$-separated set for $h^{-1}$ with cardinality $r_{n(p-1)}(\varepsilon, h^{-1})$. Fix any $x \in \tilde{\Psi} \cap [0,1]^2$ and define $E = \{(x, c) : c \in C_n\}$. We claim that $C_n$ is an $(n, \delta)$-separated set for $\tilde{H}^k_C$.

Take any two distinct points $(x, r), (x, t) \in E$ and let $0 \leq j < n(p+1)$ be such that $d(h^j(r), h^j(t)) > \varepsilon$. Since by (6.1) we have $\tilde{H}^{nk}(x) \in [0,1] \times [n(p-1), \infty)$, by definition of $s$ there is $0 \leq i < n$ such that if we denote $\tilde{H}^{ik}(x) = (q, u)$ then $u - j \in (-2s, 2s)$. Let $w \geq 0$ be an integer such that $u \in [w, w+1)$. Then
$$\tilde{H}^{ik}_C((x,r)) \approx (\tilde{H}^{ik}(x) - w, h^{-w}(t)) \in D,$$
$$\tilde{H}^{ik}_C((x,t)) \approx (\tilde{H}^{ik}(x) - w, h^{-w}(r)) \in D$$
and clearly $d(h^{-w}(t), h^{-w}(r)) > \gamma$ because otherwise $d(h^{-j}(t), h^{-j}(r)) < \varepsilon$ which is a contradiction. But then
$$\delta < d(Q(\tilde{H}^{ik}_C((x,r)), Q(\tilde{H}^{ik}_C((x,t)))) = d(H^{ik}_C(Q(x,r)), H^{ik}_C(Q(x,t)))$$
which shows that the set $Q(E)$ is $(n, \delta)$-separated for $\tilde{H}^k_C$. We have the following
$$\limsup_{n \to \infty} \frac{1}{n} \log s_n(H^k_C, \delta) \geq \limsup_{n \to \infty} \frac{1}{n} \log(s_{n(p-1)}(\varepsilon, h^{-1}))$$
$$\geq \limsup_{n \to \infty} \frac{1}{n} \log(s_n(\gamma, h^{-(p-1)})$$
which in turn gives $h_{\text{top}}(H_C) \geq \frac{p-1}{k} h_{\text{top}}(h^{-1})$ and finally we obtain desired inequality
$$h_{\text{top}}(H_C) \geq \alpha h_{\text{top}}(h^{-1})$$
which completes the proof. $\square$

**COROLLARY 6.2.** *For every $\beta \in [0, \infty]$ there exists a one-dimensional hereditarily indecomposable continuum $M$ and a transitive homeomorphism $F \colon M \to M$ with $h_{top}(F) = \beta$.*

*Proof.* If $\beta = 0$ then it is enough to use homeomorphism $F$ from example of Handel [44], and if $\beta = +\infty$ then we may refer to any homeomorphism $F$ of the pseudo-circle constructed using Minc-Transue type construction from [57] (see also [19]). Therefore, we may assume that $\beta \in (0, \infty)$.

Use method of Handel [44] to construct a homeomorphism $(\Psi, H|_\Psi)$ of the pseudo-circle $\Psi$ with a rotation number $\alpha > 0$. Take a Bernoulli automorphism defined by a measure $\nu$ with entropy $h_\nu = \beta/\alpha$. Every Bernoulli automorphism is a Kolmogorov automorphism, in particular $\mu$ is strongly mixing. Since entropy of a measure is an isomorphism invariant, by the Jewett-Krieger theorem we obtain a homeomorphism $f \colon C \to C$ of a Cantor set $C$ which is uniquely ergodic with respect to $\mu$, and so if we out $h = f^{-1}$ then $h_{\text{top}}(h^{-1}) = h_\mu(f) = \beta/\alpha$ and $h$ is minimal and mixing (e.g. see [77, §4.9] for more details). Now it is enough to apply Theorems 6.1, 3.6 and 4.1 to end the proof to see that $F = H_C$ is transitive with $h_{\text{top}}(F) = \beta$ and $\Psi_C$ is hereditarily indecomposable. $\square$

## 7. Admissible values of entropy

It is well known that the rotation set of annulus map is a topological invariant [21]. The situation is a bit more delicate for cofrontier maps, as it is not immediately clear that the rotation set is independent of embedding, in especially if the rotation set consists of a single irrational number, and although any cofrontier homeomorphism extends to an annulus map it is obvious that distinct extensions could produce distinct rotation sets. However, it follows from the arguments used in [21, P1. p. 257] that the rotation sets are indeed independent of embeddings. For completeness sake we sketch the argument from [21].

**PROPOSITION 7.1.** *Suppose $C$ and $K$ are two topologically inequivalent embeddings of a cofrontier in the annulus $\mathbb{A}$. Let $F : C \to C$, $G : K \to K$ and $H : C \to K$ be homeomorphisms such that $F \circ H = H \circ G$, with rotation sets $R(F)$ and $R(G)$. Then $R(F) = R(G)$.*

*Proof.* Choose lifts $f, g,$ and $h$ of $F, G$ and $H$ respectively. Let $f = (f_1, f_2)$, $g = (g_1, g_2)$, $h = (h_1, h_2)$. Note that the difference $h_1(f^n(\tilde{x})) - f_1^n(\tilde{x})$ is uniformly bounded for all $\tilde{x} \in C$ and there exists a $\tilde{y} \in K$ such that $\tilde{y} = \tilde{h}(\tilde{x})$. Therefore for any $n$ we have the following
$$\lim_{n \to \infty} \frac{f_1^n(\tilde{x})}{n} = \lim_{n \to \infty} \frac{f_1^n(\tilde{x}) + h_1(f^n(\tilde{x})) - f_1^n(\tilde{x})}{n} = \lim_{n \to \infty} \frac{h_1(f^n(h^{-1}(\tilde{y})))}{n}.$$
It follows that $R(F) = R(H^{-1} \circ F \circ H) = R(G)$. $\square$

Now recall that for every rational number $\frac{p}{q}$ there exists an annulus homeomorphism with an invariant pseudo-circle such that its rotation set is $\{\frac{p}{q}\}$. The same holds for its restriction to the pseudo-circle. Indeed, fix a rational $\frac{p}{q}$ and let the pseudo-circle $\Psi$ be embedded in an annulus $\mathbb{A}$. Then the $q$-fold covering $(A_q, \tau_q)$ contains the $q$-fold covering $\Psi_q = \tau_q^{-1}(\Psi)$ of $\Psi$, which is homeomorphic to the pseudo-circle [46]. Let $\sigma$ be a generator of the deck transformation group of $(A_q, \tau_q)$. Then $\sigma^p \colon A_q \to A_q$ is a periodic homeomorphism with $\sigma^p(\Psi_q) = \Psi_q$. It is easy to see, from the fact that $\sigma$ is a deck transformation, that the rotation set of $\sigma^p$ is equal to $\{\frac{p}{q}\}$, and it is the same for $\sigma^p|\Psi_q$. By Kerékjártó's Theorem [28] $\sigma^p$ is conjugate to an isometry. Consequently, without violating any essential properties of HAKs that are required in the proof of Theorem 6.1, on the pseudo-circle we can compose a HAK homeomorphism with a periodic homeomorphism with unique rotation number $\frac{p}{q}$ to get the following corollary.

**COROLLARY 7.2.** *For every $\beta \in (0, \infty)$ and every $t \in \mathbb{Q}_+$ there exists a one-dimensional hereditarily indecomposable continuum $X_\beta$ and a homeomorphism $F_t \colon X_\beta \to X_\beta$ such that $h_{top}(F_t) = t \cdot \beta$. In particular, there exists an uncountable collection of one-dimensional hereditarily indecomposable continua, each of which admits a dense set in $\mathbb{R}$ of entropy values for its homeomorphisms.*

We say that $M$ is a *Mycielski set* if it can be presented as a countable union of Cantor sets.

**THEOREM 7.3.** *There is a dense Mycielski set $M \subset \mathbb{R}$ such that each $\alpha \in M$ can be obtained as the rotation number of a minimal homeomorphism in Handel's construction.*

*Proof.* Fix any nonempty open set $U \subset \mathbb{R}$. We are going to show that there exists a Cantor set $C \subset U$ such that each number in $C$ represents the rotation number of a homeomorphism obtained by Handel's method. Fix any $\gamma \in U$ and any $\varepsilon > 0$ such that $[\gamma - \varepsilon, \gamma + \varepsilon] \subset U$. The rotation number of Handel's example $\alpha = \lim_{n \to \infty} \alpha_n$ is obtained as a sequence of approximations in such a way that $\alpha_1 = \gamma$ and consecutive numbers $b_n = \alpha_{n+1} - \alpha_n$ are sufficiently small. We can always decrease $b_n$ (this may also force us to decrease all $b_i$ for $i > n$), but it can be done inductively. Therefore we can select a sequence $b_n$ in such a way that if we denote $b_n^1 = b_n$ and $b_n^0 = b_n/3$ then

(1) each $b_n > 0$,
(2) $\sum_{n=1}^{\infty} b_n < \varepsilon/2$,
(3) $\sum_{n=k+1}^{\infty} b_k < b_n/6$,
(4) for any sequence $i_1, \ldots, i_n \in \{0,1\}$ the sequence of rotations $\alpha_n = \sum_{j=1}^{n} b_j^{i_j}$ can be used in Handel's construction in [44] (all conditions required in this construction up to step $n$ are satisfied).

By the above conditions, for any sequence $x \in \{0,1\}^\mathbb{N}$ we can construct a Handel's homeomorphism with the rotation number $\alpha_x = \lim_{n \to \infty} \sum_{k=1}^{n} b_k^{x_k}$. It is not hard to see that the assignment $x \mapsto \alpha_x$ is continuous and that if $x \neq y$ then $\alpha_x \neq \beta_x$. Therefore we obtain a closed and uncountable set of rotation numbers of Handel's homeomorphism in $U$. This set clearly contains a cantor set finishing the first part of the proof. But now, to construct $M$ it is enough to take union of all Cantor sets that can be constructed for sets from a countable basis of open set $U$ for the topology on $\mathbb{R}$. $\square$

Extending the above observations on $k$–fold coverings of HAK's, we obtain the following.

**COROLLARY 7.4.** *For every $\beta \in (0, \infty)$ there exists a dense Mycielski set $M \subset (0, +\infty)$ and a one-dimensional hereditarily indecomposable continuum $X_\beta$ such that:*

(1) $\beta \in M$ and $\mathbb{Q}_+ M = M$,
(2) for every $t \in M$ there is a weakly mixing homeomorphism $F_t \colon X_\beta \to X_\beta$ such that $h_{top}(F_t) = t$.

## 8. Smooth extensions

To form the pseudo–suspension we endowed $C$ with a $\mathbb{Z}$–space structure induced by a homeomorphism $h \colon C \to C$. The minimal homeomorphisms of $C$ used in the construction of the minimal maps of $\Psi_C$ in are subshifts of $\{0,1\}^{\mathbb{Z}}$. Pick any such subshift $h$ and consider the full shift $s$ of $\{0,1\}^{\mathbb{Z}}$. It is well known that there is a $C^\infty$–smooth diffeomorphism $f : M \to M$ of a compact, 2–dimensional manifold $M$ with an invariant Smale horseshoe $\Lambda$ with restriction $f|\Lambda$ that is topologically conjugate to $s$. Thus, there is an invariant, minimal subset $K \subset M$ such that $f|K$ is topologically conjugate to $h$.

As in Subsection 3.1, we now begin with the principal $\mathbb{Z}$–bundle $\xi = ([0,1] \times \mathbb{R}, q, \mathbb{A})$. We then give $M$ the $\mathbb{Z}$–space structure induced by $f$, $n.x \mapsto f^{-n}(x)$, giving us the fiber bundle $\xi[M]$. The total space $([0,1] \times \mathbb{R})_M$ of $\xi[M]$ is the quotient space $([0,1] \times \mathbb{R}) \times M/\approx$, where $((s,r),x) \approx ((s',r'),x')$ if and only if $s = s'$ and there is an $n \in \mathbb{Z}$ satisfying $r' = r + n$ and $x' = f^{-n}(x)$, and we denote the quotient map as $Q_M : ([0,1] \times \mathbb{R}) \times M \to ([0,1] \times \mathbb{R})_M$.

A direct way of describing the smooth structure of $([0,1] \times \mathbb{R})_M$ is in terms of induced functional structures, see, e.g, [26, II,2 and 4].

**DEFINITION 8.1.** *A* functional structure *on a topological space $X$ is a function $F_X$ on the collection of open subsets of $X$ satisfying:*

(1) *For each open $U \subset X$, $F_X(U)$ is a subalgebra of all continuous real valued functions on $U$;*
(2) *$F_X(U)$ contains all constant functions;*
(3) *For open $U \subset V$, if $f \in F_X(U)$, then $f|V \in F_X(V)$; and*
(4) *For open $U = \cup U_\alpha$ and a continuous function $f$ on $U$, if $f|U_\alpha \in F_X(U_\alpha)$ for all $\alpha$, then $f \in F_X(U)$.*

*And in this case $(X, F_X)$ is called a* functionally structured space. *If $U$ is an open subset of $X$, then $(U, F_U)$ is the functional structure of $U$ induced by $F_X$.*

*A $C^\infty$ $n$–dimensional manifold with boundary $X$ is a second countable functionally structured Hausdorff space such that for each point $x \in X$ there is a neighbourhood $U \subset X$ of $x$ and there is a corresponding open set $V \subset \mathbb{R}^n_+$ such that $(U, F_U)$ is isomorphic to $(V, C^\infty)$, where $C^\infty$ assigns to $V$ the $C^\infty$ functions on $V$. Here $\mathbb{R}^n_+$ is the set of points in $\mathbb{R}^n$ with non–negative first coordinate, and two functionally structured spaces $(Y, F_Y), (Z, F_Z)$ are* isomorphic *if there is a homeomorphism $h : Y \to Z$ such that $F_Y = F_Z \circ h$, where $h$ denotes the function $h$ induces on the collection of open subsets.*

Now consider a $C^\infty$ functional structure of $([0,1] \times \mathbb{R}) \times M$ obtained from the product structure, and endow $([0,1] \times \mathbb{R})_M$ with the functional structure induced by $Q_M$, $(([0,1] \times \mathbb{R})_M, F_{([0,1] \times \mathbb{R})_M})$. That is, for an open $U \subset ([0,1] \times \mathbb{R})_M$,

$$(g : U \to \mathbb{R}) \in F_{([0,1] \times \mathbb{R})_M}(U) \text{ iff } g \circ Q_M \in F_{([0,1] \times \mathbb{R})}\left(Q_M^{-1}(U)\right).$$

This will give $([0,1] \times \mathbb{R})_M$ a functional structure that is locally isomorphic to that of $\mathbb{R}^4_+$ since two sufficiently small open sets in $([0,1] \times \mathbb{R}) \times M$ that are mapped to the same open set in $([0,1] \times \mathbb{R})_M$ will differ by translations of the $\mathbb{R}$ factor and a corresponding iterate of the $C^\infty$ map $f$ in the $M$ factor, see, [26, II,Ex. 4.2]. Thus, $([0,1] \times \mathbb{R})_M$ is a $C^\infty$ manifold of dimension 4 with boundary with this functional structure. Just as in Subsection 3.2 we have the fiber bundle projection $q_M : ([0,1] \times \mathbb{R})_M \to \mathbb{A}$ which has the homotopy lifting property. We give $\mathbb{A}$ the $C^\infty$ structure induced by $q_M$ and construct a Handel diffeomorphism $H : \mathbb{A} \to \mathbb{A}$ with an attracting pseudo-circle $\Psi$ as before. We then have the lifted homeomorphisms $H_M, \widetilde{H_M}$ as in Equations 3.1 and 3.2

(8.1)
$$\begin{array}{ccc} [0,1]\times\mathbb{R}\times M & \xrightarrow{\widetilde{H_M}} & [0,1]\times\mathbb{R}\times M \\ Q_M \downarrow & & \downarrow Q_M \\ ([0,1]\times\mathbb{R})_M & \xrightarrow{H_M} & ([0,1]\times\mathbb{R})_M \\ q_M \downarrow & & \downarrow q_M \\ \mathbb{A} & \xrightarrow{H} & \mathbb{A}. \end{array}$$

Notice that $\widetilde{H_M}$ is $C^\infty$–smooth since locally it is the product of the lift of $H$ to the universal cover $[0,1]\times\mathbb{R}$ of $\mathbb{A}$ and the identity on $M$. As we have endowed $([0,1]\times\mathbb{R})_M$ with the smooth structure induced by $Q_M$, $H_M$ is then equally smooth. Now $H_M$ has a topological copy of $\Psi_C$ as an invariant minimal set with intermediate complexity, and by choosing the rotation number of $H$ as described in Section 6, the restricted map on $\Psi_C$ can have entropy of all values as described in Corollary 7.2.

Piecing this together, we finally obtain the main theorem.

**THEOREM 8.2.** *There exist hereditarily indecomposable continua admitting homeomorphisms of intermediate complexity. Moreover, there exist an hereditarily indecomposable continuum $\Psi_C$ satisfying:*

(1) *$\Psi_C$ occurs as an invariant minimal set in a smooth diffeomorphism $F$ of the 4 dimensional manifold $([0,1]\times\mathbb{R})_M$ and*
(2) *The restriction $F|\Psi_C$ is weakly mixing homeomorphism of intermediate complexity.*

**REMARK 8.3.** *In the above theorem, it should be observed that for a* single *such diffeomorphism $F$ as above, we obtain any continuum that can be generated by a subshift of $\{0,1\}^\mathbb{Z}$ since this shift is conjugate to a subset of Smale's horseshoe.*

**REMARK 8.4.** *If the diffeomorphism $f\colon M\to M$ is isotopic to the identity, then the manifold $([0,1]\times\mathbb{R})_M$ will be diffeomorphic to $M\times\mathbb{A}$. In particular, if we start with $M$ as a disc, then we can easily obtain $\Psi_C$ as an invariant set of a diffeomorphism on $\mathbb{R}^4$.*

In the above we were aiming to find smooth maps for which $\Psi_C$ occurs as an invariant set within an attractor, which is based on the Handel example in which $\Psi$ occurs as an attractor. But Handel also constructs examples for which $H\colon\mathbb{A}\to\mathbb{A}$ preserves area. It is such examples that we now use to construct another family of examples.

*Proof of Theorem 1.2.* For this example we take a hyperbolic automorphism $f:\mathbb{T}^2\to\mathbb{T}^2$ to construct our diffeomorphism of $([0,1]\times\mathbb{R})_{\mathbb{T}^2}$ as above. Using Markov partitions, one can construct a subshift of finite type $\Omega_A$ and a finite–to–one factor map $g:\Omega_A\to\mathbb{T}^2$ (see e.g. Theorems 4.3.5. and 4.3.6 in [2].) Since $\Omega_A$ is a mixing shift of finite type with positive entropy, using standard methods (e.g. the one introduced by Grillenberger in [43, Section §2]), we can find a mixing strictly ergodic (hence minimal) subshift $W\subset\Omega_A$ with positive entropy. The image $g(W)$ will be a minimal, invariant subset of $\mathbb{T}^2$ and so clearly also weakly mixing. But $g$ is finite-to-one, and hence the entropy of $f$ on $g(W)$ is the same as $\sigma$ on $W$ (first proved by Bowen [22], but see also [33, Theorem 7.1]). As before, this yields the desired hereditarily indecomposable continuum $\Psi_C$ with associated weakly mixing homeomorphism of intermediate complexity and we need only show that the diffeomorphism $H_{\mathbb{T}^2}$ preserves volume. But $H$ preserves the area on $\mathbb{A}$ and $f$ preserves the area element of $\mathbb{T}^2$ that is induced from the area element of $\mathbb{R}^2$ by the covering map $\mathbb{R}^2\to\mathbb{T}^2$. Thus, if we endow $([0,1]\times\mathbb{R})_{\mathbb{T}^2}$ with the volume element induced by the quotient map $Q_M$ using the indicated area elements in the product covering space, we see that $H_M$ preserves volume. □

Observe that in the above example $f$ is not homotopic to the identity, making it more difficult to describe the manifold $([0,1]\times\mathbb{R})_{\mathbb{T}^2}$ in more basic terms. Also, there will be a plethora of subshifts $W$ that one can take, leading to many different invariant hereditarily indecomposable continua within $([0,1] \times \mathbb{R})_{\mathbb{T}^2}$ with the properties listed above. However, it will be a more restrictive collection since not all subshifts would be included as in the previous construction.

It is also worth mentioning that the Cantor set can never appear as an attractor in a manifold (or even locally connected compact metric space).[1] However many subshifts can be obtained as isolated invariant sets (i.e. only sets satisfying $f(A) = A$ in some open neighborhood of $A$) for diffeomorphisms on manifolds [25].

## 9. Pseudo-solenoids

In this section we shall investigate the lift of Handel's homeomorphism to pseudo–solenoids and establish some of its properties: minimality, weak mixing, uniform rigidity and zero entropy. Recall that a pseudo–solenoid can be viewed as a principal bundle over the pseudo–circle with the bundle structure that is induced by the principal bundle structure of a standard Vietoris solenoid over a circle, see [69]. We shall now put this bundle structure into the context of the pseudo–suspension as developed earlier. For a given sequence $P = (p_1, p_2, p_3, \ldots)$ of prime numbers $p_i$, there is the corresponding adding machine $a_P : C \to C$ and the corresponding solenoid

$$\Sigma_P = \varprojlim \{(\mathbb{R}/\mathbb{Z}, p_i) : i \in \mathbb{N}\}$$

obtained by taking the inverse limit of the covering maps the circle $[x] \xmapsto{p_i} [p_i x]$. Moreover, the solenoid $\Sigma_P$ can be viewed as the suspension of the matching adding machine $a_P$, [1]. The $P$-adic pseudo-solenoid $\Psi_P$ is then the total space of the induced bundle as indicated in the below Diagram that corresponds to the earlier Diagram 3.3

(9.1)
$$\begin{array}{ccc} \Psi_P & \xrightarrow{F} & \Sigma_P \\ {\scriptstyle q_P}\downarrow & & \downarrow{\scriptstyle q_{a_P}} \\ \Psi & \xrightarrow{f} & \mathbb{R}/\mathbb{Z} \end{array}$$

and as such is the pseudo–suspension of $a_P$. But as the bundle for $\Sigma_P$ is principal, the same will be true for the $\Psi_P$ bundle map, unlike the typical general fiber bundle in a pseudo–suspension. And while we have the lifted homeomorphism $H_P : \Psi_P \to \Psi_P$ as before, because the adding machine $a_P$ is equicontinuous, not weakly mixing and not totally minimal, none of the arguments we previously used to establish the weak mixing and minimality of the lifted homeomorphism will apply to this case. Instead, we will have to exploit the alternative representation of $\Psi_P$ as an inverse limit

$$\Psi_P = \varprojlim \{(\Psi, \tau_i) : i \in \mathbb{N}\}$$

where $\tau_i$ is a covering map of $\Psi$ of degree $p_i$. Consideration of the lifting construction of $H_P$ and the nature of the inverse limit, leads to the following commutative diagram

(9.2)
$$\begin{array}{ccccccccccc} \Psi_P: & & \Psi & \xleftarrow{\tau_1} & \Psi & \xleftarrow{\tau_2} & \cdots \longleftarrow & \Psi & \xleftarrow{\tau_i} & \Psi & \longleftarrow \cdots \\ & & {\scriptstyle H_P}\downarrow & & {\scriptstyle H}\downarrow & & {\scriptstyle H_1}\downarrow & & {\scriptstyle H_{i-1}}\downarrow & & {\scriptstyle H_i}\downarrow \\ \Psi_P: & & \Psi & \xleftarrow{\tau_1} & \Psi & \xleftarrow{\tau_2} & \cdots \longleftarrow & \Psi & \xleftarrow{\tau_i} & \Psi & \longleftarrow \cdots \end{array}$$

where each $H_i$ is the lift of $H$ via the covering map $\tau_i \circ \cdots \circ \tau_1$ obtained as before using the lifting property of fibrations. Thus, $H_P$ is the inverse limit of the homeomorphisms indicated in the Diagram. Each covering map $\tau_i \circ \cdots \circ \tau_1$ extends to coverings of $\mathbb{A}$ of the same degree, and so each of the $H_i$ is a Handel homeomorphism that can be realized using the intersections of crookedly nested annuli obtained from lifting the nested annuli in the base copy of $\mathbb{A}$. As each of the lifted maps $H_i$ is again a Handel homeomorphism, it shares all the qualitative properties of the original $H$, though the rotation number may be altered.

---

[1] We are grateful to Mike Boyle for bringing this this fact to our attention.

**PROPOSITION 9.1.** *Let $X = \varprojlim\{(X_i, f_i): i \in \mathbb{N}\}$ with each $X_i$ a compact metric space and all bonding maps $f_i: X_{i+1} \to X_i$ continuous surjections. Suppose that $h: X \to X$ is a homeomorphism given by $h((x_i)) = (h_i(x_i))$, where each $h_i: X_i \to X_i$ is a homeomorphism. Then,*

(1) *$h$ is weakly mixing if and only if each $h_i$ is weakly mixing,*
(2) *$h$ is minimal if and only if each $h_i$ is minimal and*
(3) *$h$ is uniformly rigid if and only if each $h_i$ is uniformly rigid.*

*In particular, for each $P$, $H_P$ is weakly mixing, minimal and uniformly rigid.*

*Proof.* The proofs of (1) and (2) are standard and can be found, for example, in [34, IV,Prop 1.5, p.275] and so we proceed to (3). Assume then that each $h_i$ is uniformly rigid. Using standard arguments of metrics for inverse limits, one can find for each $\varepsilon > 0$ a corresponding $N$ and $\varepsilon_N > 0$ so that if the $N$–th coordinates of two points of $X$ differ by no more than $\varepsilon_N$ in $X_N$, then the two points differ by no more than $\varepsilon$ in $X$. Let $\varepsilon > 0$ then be given. As the homeomorphism $h_N$ is uniformly rigid, there exists an (arbitrarily large) $k$ so that $h_N^k$ is within $\varepsilon_N$ of $\text{id}_{X_N}$ in the uniform metric. Hence, $h^k$ is within $\varepsilon$ of $\text{id}_X$. Thus, $h$ is uniformly rigid. It is clear that the uniform rigidity of $h$ implies that of each $h_i$. Since we have established each of the three properties for the general Handel homeomorphism, the result follows. □

Thus, we now have a large family of minimal, uniformly rigid and weakly mixing homeomorphisms found quite easily using inverse limit representations. It is natural to wonder whether similar techniques could be applied to the more general examples we considered earlier to establish these same properties. In [38] Fearnley proved that any two hereditarily indecomposable circularly chainable continua are homeomorphic iff their one dimensional Čech cohomology groups are isomorphic. Moreover, the family of circle–like hereditarily indecomposable continua is characterized as the family of all pseudo-solenoids $\Psi_P$ as described here, the psuedo–circle and the pseudo–arc (cf. results of Fearley [39] and Rogers, Jr. [71]). The pseudo–circle is circle–like and hence the same is true of any inverse limit of pseudo–circles. As any pseudo–suspension will have non–trivial cohomology, we know that the only continua other than the pseudo–circle that we will obtain in an inverse limit of Handel homeomorphisms of pseudo–circles are precisely these pseudo–solenoids. In a similar vein, while this inverse limit approach could be applied to finite sheeted covering spaces of pseudo–solenoids, this would not result in any new examples as any finite-sheeted covering space of $\Psi_P$ is again homeomorphic to $\Psi_P$ [20].

**THEOREM 9.2.** *For each $P$, the homeomorphism $H_P$ of the $P$-adic pseudo-solenoid $\Psi_P$ is uniformly rigid, minimal and weakly mixing. Moreover, there is a volume preserving diffeomorphism $F$ of a 4–dimensional manifold which has an invariant set $X$ that is homeomorphic to $\Psi_P$ and $F|X$ is topologically conjugate to $H_P$.*

*Proof.* In light of Proposition 9.1, we need only find the smooth example as stated. To construct the diffeomorphism $F$ as in the statement of the theorem, we apply a construction very similar to that in Theorem 1.2, only instead of a hyperbolic map of a torus, we choose an area preserving diffeomorphism of a disk with an invariant adding machine $a_P$ that can be extended to $S^2$, (e.g. see [30]). □

The diffeomorphism of $S^2$ in Theorem 9.2 can be constructed to be isotopic to the $id_{S^2}$, and so the manifold can be chosen to be $S^2 \times \mathbb{A}$. The fact that $H_P$ has zero entropy follows directly from the fact that it is uniformly rigid, hence does not have asymptotic pairs (cf. [14]).

## 10. Nondegenerate rotations sets and infinite entropy - Proof of Theorem 1.3

Let $\mathcal{U} = [U_1, ..., U_7]$ be a taut chain cover of a continuum $X$. Following [63] we say that $H \colon X \to X$ *stretches* $\mathcal{U}$ if $H(U_3) \subset U_1$ and $H(U_5) \subset U_7$ (or $H(U_5) \subset U_1$ and $H(U_3) \subset U_7$). We say that $H$ is a *stretching map* of $X$ if there exists a taut chain cover $\mathcal{U}$ of $X$ and positive integer $n$ such that $H^n$ stretches $\mathcal{U}$.

A function $f \colon \{1, \ldots, m\} \to \{1, \ldots, n\}$ is called a *pattern* provided $|f(i+1) f(i)| \leq 1$ for $i \in \{1, \ldots, m-1\}$. A pattern $f \colon \{1, \ldots, 2k+5\} \to \{1, \ldots, 7\}$ is a *k-fold* provided that $k$ is odd and

$$f(i) = \begin{cases} i, & i \leq 5 \\ 3, & i = 3 + 4m \text{ for some } m \in \mathbb{N}, \\ 4, & i = 4 + 2m \text{ for some } m \in \mathbb{N}, \\ 5, & i = 5 + 4m \text{ for some } m \in \mathbb{N}, \\ 6, & i = 2k + 4, \\ 7, & i = 2k + 5. \end{cases}$$

The following Theorem is a very useful tool, proved first by Oversteegen and Tymchatyn in [66]:

**THEOREM 10.1.** *Let $C$ be a hereditarily indecomposable continuum and let $\mathcal{U} = \{U_1, \ldots, U_n\}$ be an open taut chain cover of $C$. Let $f \colon \{1, \ldots, m\} \to \{1, \ldots, n\}$ be a pattern on $\mathcal{U}$. Then there exists an open taut chain cover $\mathcal{V}$ of $C$ such that $\mathcal{V}$ follows the pattern $f$ in $\mathcal{U}$.*

The following fact is a simple extension of the result in [63] inspired by an earlier paper by Kennedy and Yorke [55]. The idea of the proof is essentially the same.

**THEOREM 10.2.** *Let $U = \{U_1, \ldots, U_7\}$ be a taut chain cover of the pseudo-arc $C$ and assume that $f \colon C \to C$ stretches $\mathcal{U}$. Then, for every $k \geq 2$ there exists an invariant set $D \subset \overline{U_4}$ such that $(D, f|_D)$ is an extension of the full shift on $k$ symbols $(\Sigma_k^+, \sigma)$.*

*Proof.* Fix any odd integer $k \geq 2$. Let $\mathcal{V} = \{V_1, \ldots, V_{2k+5}\}$ be a $k$-fold refinement of $\mathcal{U}$. For every index $i_0 \in \{1, \ldots, k\}$ there exists a continuum $Y_{i_0} \subset V_{2i_0+2}$ such that $Y_{i_0} \cap V_{2i_0+1} \neq \emptyset$ and $Y_{i_0} \cap V_{2i_0+3} \neq \emptyset$.

Since $H(Y_{i_0})$ is a continuum intersecting both $V_1$ and $V_{2k+5}$, for each index $i_0 \in \{1, \ldots, k\}$ there exists a continuum $Y_{i_0, i_1} \subset H(Y_{i_0}) \cap V_{2i_1+2}$ such that $Y_{i_1} \cap V_{2i_1+1} \neq \emptyset$ and $Y_{i_1} \cap V_{2i_1+3} \neq \emptyset$.

Proceeding inductively, for every finite sequence of indexes $i_0, i_1, \ldots, i_n$ there exist continua $Y_{i_0}, \ldots, Y_{i_0, \ldots, i_n}$ such that $Y_{i_0, \ldots, i_{k+1}} \subset H(Y_{i_0, \ldots, i_k}) \cap V_{2i_k+1}$ for each $k = 0, \ldots, n-1$. This implies that

$$\bigcap_{k=0}^{n} H^{-k}(Y_{i_0, \ldots, i_k}) \neq \emptyset$$

and so there exists a point $z_{i_0, \ldots, i_n}$ such that we have $H^k(z_{i_0, \ldots, i_n}) \in Y_{i_0, \ldots, i_k}$ for each $k = 0, \ldots, n$. Using compactness, for every $x \in \Sigma_k^+$ we can find a point $z_x$ such that $H^k(z_x) \in \overline{V_{2x_k+2}}$. For every $x \in \Sigma_k^+$ denote by $D_x$ the set such that $H^k(z) \in \overline{V_{2x_k+2}}$ for every $z \in D_x$ and every $k \geq 0$. Clearly each $D_x$ is closed and nonempty. Furthermore, directly from definition we have that $H(D_x) \subset D_{\sigma(x)}$ and that the sets $D_x$ are pairwise disjoint.

Denote $D = \bigcup_{x \in \Sigma_k^+} D_x$. Note that for every $k$ if $\lim_{j \to \infty} p_j = q$ and each $H^k(p_j) \in \overline{V_{2i+2}}$ for some $i$ then also $H^k(q) \in \overline{V_{2i+2}}$. This shows that if $z \in \overline{D}$ then there is always $x \in \Sigma_k^+$ such that $z \in D_x$. Therefore $D = \overline{D}$. Define $\pi \colon D \to \Sigma_k^+$ by putting $\pi(z) = x$ whenever $z \in D_x$. By the previous observations we have that $\pi$ is continuous, surjective, and $\pi \circ H = \sigma \circ \pi$. The proof is completed by the fact that $D \subset \bigcup_{i=1}^k \overline{V_{2i+2}} \subset \overline{U_4}$. $\square$

Recall the universal covering of $\mathbb{A} = [0, 1] \times \mathbb{R}/\mathbb{Z}$ is given by $q \colon [0, 1] \times \mathbb{R} \to \mathbb{A}$, where $q = \text{id}_{[0,1]} \times p$ and $p$ is the quotient map $\mathbb{R} \to \mathbb{R}/\mathbb{Z}$. A continuum $E \subset \text{Int}\,\mathbb{A}$ is *essential* if $\mathbb{A} \setminus E$ has two connected components $\mathcal{U}_+, \mathcal{U}_-$, each of which contains a different component of the boundary of $\mathbb{A}$.

**LEMMA 10.3.** *Let $A \subset \mathbb{A}$ be an essential annulus containing the pseudo-circle $\Psi$ and let $H\colon \mathbb{A} \to \mathbb{A}$ be a homeomorphism such that $H(A) \subset \operatorname{Int} A$. Assume that there are periodic points $p_0, p_1 \subset A$ with different rotation numbers and two contractible continua $C_0, C_1 \subset \mathbb{A}$ such that:*

*(1) $p_i \in C_i$ and $H^{n_i}(C_i) \subset C_i$ where $n_i$ is the period of $p_i$,*
*(2) $C_i$ is inessential and intersects both boundary components of $A$.*

*Then $H$ has infinite entropy.*

*Proof.* Denote $\tilde{A} = [0,1] \times \mathbb{R}$. Let $g$ be an iterate of $H$ such that $p_0, p_1$ are fixed points of $g$ and let $G$ be the lift of $g$ such that $G(q_0) = q_0 - j$ and $G(q_1) = q_1 + k$ where $j, k$ are positive integers and $q_i$ is a lift of $p_i$.

There exists a continuum $D_i \subset C_i$ which intersects both boundary components of $A$ (see Theorem 14.3 in [65]). Note that $G(\tilde{C}_0) \subset \tilde{C}_0 - j$ and $G(\tilde{C}_1) \subset \tilde{C}_1 + k$, here $\tilde{C}_i$ is a lift of $C_i$. Then, if we denote by $\tilde{D}_i$ a lift of $D_i$ then there exists $m > 0$ such that $G^m(\tilde{D}_0)$ is on the left of $\tilde{D}_0, \tilde{D}_1$ and $G^m(\tilde{D}_1)$ are on the right of $\tilde{D}_0, \tilde{D}_1$. Let $\tilde{\Psi}$ be a lift of $\Psi$ and denote by $P$ its two point compactification by points $-\infty, \infty$ which is homeomorphic to the pseudo-arc (see [58], [9]). Then $G$ extends to $P$ by putting $G(-\infty) = -\infty$ and $G(\infty) = \infty$. Denote by $U_4$ the intersection of $P$ with the open set contained between the continua $\tilde{D}_0, \tilde{D}_1$ in $A$ and denote by $U_3$ and $U_5$ the intersections of $P$ with sufficiently small neighborhoods of $\tilde{D}_0, \tilde{D}_1$ respectively, such that there are neighborhood $U_1$ of $-\infty$ and $U_7$ of $\infty$ satisfying $G^m(U_3) \subset U_1$ and $G^m(U_5) \subset U_7$ and $\overline{U_1} \cap \overline{U_3 \cup U_4 \cup U_5 \cup U_7} = \emptyset$ and $\overline{U_7} \cap \overline{U_1 \cup U_3 \cup U_4 \cup U_5} = \emptyset$. Then there are open sets $U_2, U_6$ such that $U_1, \ldots, U_7$ is a taut chain cover of $P$. This shows that $G$ is a stretching map of $P$. By Theorem 10.2, for every $l \geq 2$ we obtain a $G^m$-invariant set $Q_l \subset P$ such that $G^m$ on $Q_l$ is an extension of full shift on $l$ symbols. Clearly we may assume that $Q_l \subset \tilde{\Psi}$ for each $l$, since we can always remove the fixed points $-\infty, \infty$ from $Q_l$ by removing at most two symbols from the alphabet of the corresponding shift $\Sigma_l$. Since $Q_l$ is bounded in $\tilde{\Psi}$, the map $\pi|_{Q_l}$ is finite-to-one and hence $h_{\operatorname{top}}(H) \geq h_{\operatorname{top}}(H|_{Q_l}) = \frac{1}{m}\log(l)$, which completes the proof. $\square$

The following fact is special case of Proposition 3.9 in [67].

**LEMMA 10.4.** *Let $\Psi \subset \mathbb{A}$ be the pseudo-circle attracting all the points from $\operatorname{Int} \mathbb{A}$ and assume that $H\colon \mathbb{A} \to \mathbb{A}$ is a homeomorphism with a nondegenerate rotation set. Then there are periodic points $p_0, p_1 \in \Psi$ and inessential contractible continua $C_0, C_1$ such that $p_i \in C_i$, $H^{n_i}(C_i) \subset C_i$ where $n_i$ is the period of $C_i$, $C_0 \cap C_1 = \emptyset$ and $C_i \cap \mathcal{U}_+ \neq \emptyset$, $C_i \cap \mathcal{U}_- \neq \emptyset$ for $i = 0, 1$.*

*Proof of Theorem 1.3.* Observe that statement of Lemma 10.4 complements the assumptions of Lemma 10.3, therefore we have just completed the proof of Theorem 1.3. $\square$

## 11. Concluding remarks

In [41] surprising examples are found of homeomorphisms that are minimal, uniformly rigid and weakly mixing, and our results here give us large new classes of examples satisfying these same three properties. However, it is not yet entirely clear when our lifted homeomorphisms $H_C$ satisfy these conditions. By consideration of the pseudo–solenoid, we see that the condition that $h\colon C \to C$ be weakly mixing is not necessary for the weak mixing of the lifted homeomorphism to $\Psi_C$. Similarly, since adding machines are very far from totally transitive, we see that the total transitivity of $h$ is not a necessary condition for the minimality of the lifted homeomorphism to $\Psi_C$. This naturally leads to the following question.

**QUESTION 11.1.** *Can one find necessary and sufficient conditions for $h\colon C \to C$ to guarantee that the lift of a given HAK homeomorphism $H_C$ is minimal? What about weakly mixing if the original HAK is weakly mixing?*

In our proof of minimality for the lifted homeomorphism, we used HAK homeomorphisms with rotation numbers that are very well approximated by rational numbers.

**QUESTION 11.2.** *Are there weakly mixing minimal homeomorphisms $h\colon C \to C$ such that HAK homeomorphisms with certain rotation numbers do not lift to minimal homeomorphisms $H_C$?*

Proof of Corollary 7.2 strongly relies on the possible values of rotations numbers in HAK homeomorphism. From the construction it is not clear if all values are admissible. This leads to the following question, closely related to the original question of Barge (see [61]).

**QUESTION 11.3.** *Is there a hereditarily indecomposable continuum $X$ such that for every $t \geq 0$ there is a homomorphism $F_t\colon X \to X$ of entropy $h_{top}(F_t) = t$?*

One way to address the above question in the affirmative would be to answer the following one.

**QUESTION 11.4.** *For every $\alpha \in \mathbb{R}$ is there an indecomposable cofrontier $X_\alpha$ that admits a homeomorphism $F_\alpha : X_\alpha \to X_\alpha$ with a well defined rotation number $\alpha$? Can $X_\alpha$ be realized as the pseudo-circle?*

It is also unknown wheather a minimal homeomorphism of an indecomposable cofrontier must be weakly mixing.

**QUESTION 11.5.** *(cf. [51]) Is there an indecomposable cofrontier that admits a minimal homeomorphism semi-conjugate to an irrational circle rotation?*

Note that in the proof of Theorem 1.3 we strongly rely on the fact that $\Psi$ is an attractor. It is not clear whether this is a necessity or only a technical assumption needed in our argument. This leads to the following natural question.

**QUESTION 11.6.** *Suppose that $H\colon \Psi \to \Psi$ is a homeomorphism of the pseudo-circle embedded in the plane with a nondegenerate rotation set. Is $h_{top}(H) = +\infty$?*

The answer to the above question can be an important step towards the characterization of values of entropy admissible for homeomorphism of pseudo-circle.

Finally, it would be of interest to know more about embedability of pseudo-suspensions.

**QUESTION 11.7.** *Let $h : C \to C$ be a Cantor set homeomorphism. Suppose that the suspension of $h$ embeds in a manifold $M$. Does every pseudosuspension of $h$ embed in $M$?*

## Acknowledgements

The authors are grateful to Marcy Barge and Tobias Jäger for helpful conversations. The authors are grateful to Joe Auslander, Eli Glasner and Tomasz Downarowicz for interesting discussions on construction of dynamical systems with specified properties and to Mike Boyle for his remarks on attractors and isolated sets in manifolds.

This research was supported by the project LQ1602 IT4Innovations excellence in science, Grant IN-2013-045 form the Leverhulme Trust for an International Network, which supported research visits of the authors. J. Boroński's work was supported by National Science Centre, Poland (NCN), grant no. 2015/19/D/ST1/01184 and Piotr Oprocha was supported by National Science Centre, Poland (NCN), grant no. 2015/17/B/ST1/01259.

(J. P. Boroński) AGH University of Science and Technology, Faculty of Applied Mathematics, al. Mickiewicza 30, 30-059 Kraków, Poland – and – National Supercomputing Centre IT4Innovations, Division of the University of Ostrava, Institute for Research and Applications of Fuzzy Modeling, 30. dubna 22, 70103 Ostrava, Czech Republic

*E-mail address*: `jan.boronski@osu.cz`

(A. Clark) Department of Mathematics, University of Leicester, University Road, Leicester LE1 7RH, United Kingdom

*E-mail address*: `Alex.Clark@le.ac.uk`

(P. Oprocha) AGH University of Science and Technology, Faculty of Applied Mathematics, al. Mickiewicza 30, 30-059 Kraków, Poland – and – National Supercomputing Centre IT4Innovations, Division of the University of Ostrava, Institute for Research and Applications of Fuzzy Modeling, 30. dubna 22, 70103 Ostrava, Czech Republic

*E-mail address*: `oprocha@agh.edu.pl`